\documentclass[12pt]{amsart}

\usepackage{amssymb}

\makeatletter
\@addtoreset{equation}{section}
\makeatother

\theoremstyle{plain}
\newtheorem{theorem}{Theorem}[section]
\newtheorem{lemma}[theorem]{Lemma}
\newtheorem{proposition}[theorem]{Proposition}
\newtheorem{corollary}[theorem]{Corollary}

\theoremstyle{definition}

\newtheorem{remark}[theorem]{Remark}

\newcommand{\Cee}{{\mathbb C}}
\newcommand{\Ree}{{\mathbb R}}
\newcommand{\En}{{\mathbb N}}
\newcommand{\Que}{{\mathbb Q}}
\newcommand{\Tee}{{\mathbb T}}
\newcommand{\Zee}{{\mathbb Z}}

\newcommand{\alp}{\alpha}

\newcommand{\Gam}{\Gamma}
\newcommand{\lam}{\lambda}

\newcommand{\ome}{\omega}
\newcommand{\sig}{\sigma}
\newcommand{\Sig}{\Sigma}
\newcommand{\vphi}{\varphi}

\newcommand{\fK}{\mathcal{K}}
\newcommand{\fL}{\mathcal{L}}

\newcommand{\fV}{\mathcal{V}}
\newcommand{\fW}{\mathcal{W}}

\newcommand{\rank}{\mathrm{rank}}
\newcommand{\spn}{\mathrm{span}}
\newcommand{\conv}{\mathrm{conv}}
\newcommand{\tor}{\mathrm{T}}
\newcommand{\bl}{\mathrm{L}}

\newcommand{\wbar}[1]{\overline{#1}}
\newcommand{\what}[1]{\widehat{#1}}
\newcommand{\til}[1]{\tilde{#1}}

\begin{document}

 \title[Spectra of Beurling algebras]{Spectra of Beurling algebras 
 of abelian groups}
\author{Nico Spronk}
\author{Zhihao Zhang}
\address{Dept. of Pure Mathematics\\
University of Waterloo\\
Waterloo, Ont.,~Canada\\
N2L 3G1}
\email{nico.spronk@uwaterloo.ca}
\email{z458zhang@uwaterloo.ca}

\maketitle

\begin{abstract}
We give a complete classification of all sets which are spectra of Beurling algebras
of any locally compact abelian group, parametrizing them by compact convex sets
in the dual of a certain locally convex vector space generated by the group.
 Hence, for example, a finitely generated free abelian group
may enjoy a more complicated range of such spectra than the does rational numbers.
\end{abstract}

\subsection{Introduction}\label{sec:intro}
Throughout this article we let $G$ be a locally compact abelian group. 
We use additive notation on $G$.  We consider the set of (submultiplicative,
continuous) {\em  weights}, given by
\[
W(G)=\left\{\ome:G\to(0,\infty)\,\left| \begin{matrix} \ome\text{ is continuous and} \\
\ome(x+y)\leq\ome(x)\ome(y) \\ \text{for }x,y\text{ in }G\end{matrix}\right\}\right. .
\]

\begin{remark}\label{rem:borel}
Most authors assume that a weight $\ome$ is Borel measurable, rather than continuous.
We say that two measurable weights $\ome$ and $\ome'$ are {\em equivalent} if there are $c,C>0$ so
$c\leq \frac{\ome}{\ome'}\leq C$.  Equivalent weights give rise to isomorphic Beurling
algebras, defined below.  Any weight is equivalent to a continuous one, thanks to
\cite[Theorem 3.1]{kuznetsova} (see also \cite{feichtinger}).
\end{remark}

We shall often consider the sets of weights which are {\em bounded below} or {\em symmetric},
given by
\begin{align*}
W_0(G)&=\left\{\ome\in W(G):\inf_{x\in G}\ome(x)>0\right\}\text{ and} \\
W_s(G)&=\{\ome\in W(G):\check{\ome}=\ome\}
\end{align*}
where $\check{\ome}(x)=\ome(-x)$.  

For any $\ome$ in $W(G)$ define the {\em Beurling algebra} by
\[
\bl^1(G,\ome)=\{f:G\to\Cee\;|\; f\ome\text{ is Haar integrable}\}/\sim
\]
where $\sim$ is the relation of a.e.-equivalence.  This is well-known to be a commutative
Banach algebra
under convolution, $f\ast g(x)=\int_G f(y)g(x-y)\,dy$ (Haar integral), and with norm
\[
\|f\|_{1,\ome}=\|f\ome\|_1=\int_G |f|\ome.
\]
As we shall recall below (see Theorem \ref{theo:BAspec}), the Gelfand spectrum
 of $\Sig(\bl^1(G,\ome))$, is isomorphic to a subset $\what{G}_\ome$ of the space
of complex characters $\what{G}_\Cee$.  Coincidence of this spectrum
induces an ``spectral" equivalence relation on weights.  This is rather coarse as shown by the
Gelfand-Raikov-Shilov weights in Corollary \ref{cor:GRS_minimal_weights}, below.

We will show that each weight is spectrally
 equivalent to weight of the form $\exp\circ s\circ\eta$ where $\eta:G\to\fV$
is a continuous additive map into a real locally convex
vector space, and $s$ is a  sublinear functional on $\fV$.
We will then characterize which subsets of $\what{G}_\Cee$ arise as spectra
of Beurling algebras.  All of this is summarized in Theorem \ref{theo:BAspec}.

\subsection{Organization}
Section \ref{sec:rvs} develops the theory of the universal real locally convex space
$\Ree(G)$ generated by $G$, and discusses the theory of its sublinear functionals
which are in bijective correspondence with non-empty compact convex subsets of its dual
$\Ree(G)^*$.  As appropriate, many results are developed on a general real vector space
with it finest locally convex topology.

Section \ref{sec:spectrum} relates $\Ree(G)^*$ to real characters on $G$.
Then the scope of all possible spectra of Beurling algebras $\bl^1(G,\ome)$,
for $\ome$ in $W(G)$ is resolved in Theorem \ref{theo:BAspec}.

Section \ref{sec:examples} illustrates with examples.  Most striking is Remark
\ref{rem:freab_que} in which it is indicated that $\Zee^2$ admits a more complicated
range of spectra of Beurling algebras than does $\Que$. Also of note is that 
in Remark \ref{rem:lie} we catch
a minor error which occurs in the last theorem of \cite{CartwrightMc}.

\section{Real vector spaces}\label{sec:rvs}

In this section we will discuss sets in the dual of a real vector space, and
then create the realization $\Ree(G)$ of $G$.  Certain subsets of the dual
of $\Ree(G)^*$ play a key role in our main result.

For the remainder of this section we let $\fV$  denote a real vector space.  

\subsection{The finest locally convex topology}\label{ssec:flctop}
We consider on $\fV$ the finest locally convex topology, also known as the inductive  
topology, for the collection of all finite-dimensional subspaces.  
This is described by the following equivalent descriptions:
\begin{itemize}
\item a neighbourhood basis around $0$ is given by convex, symmetric, absorbing
sets $U$ of $\fV$, for which $U\cap F$ is open in $F$ for each finite dimensional
subspace $F\cong\Ree^{\dim F}$; or
\item a neighbourhood basis around $0$ is given by $s^{-1}([0,1))$ for
seminorms $s$ on $\fV$.
\end{itemize}
Recall that given $U$ as above, the Minkowski functional 
$s_U(w)=\inf\{r>0:w\in rU\}$ is a seminorm with
$U=s_U^{-1}([0,1))$; and seminorms are automatically continuous on finite
dimensional (sub)spaces.  In particular, if $\dim\fV<\infty$,
then we take the usual topology of a finite dimensional real vector space.

\begin{remark} \label{rem:indtop}
This topology admits the following facts.
\begin{enumerate}
\item We have for a net $(v_i)$ and $v$ in $\fV$ that
\[
\lim_iv_i=v\quad\Leftrightarrow\quad \lim_i s(v_i-v)=0\quad\begin{matrix}\text{for each seminorm} \\
s\text{ on }\fV.\end{matrix}
\]
\item Hence if $\fW$ is a real locally convex vector space and $T:\fV\to\fW$ is linear,
then $T$ is continuous. In particular the space $\fV^*$ of linear functionals
consists of continuous elements.
\end{enumerate}
\end{remark}

A functional $s:\fV\to\Ree$ is called {\em sublinear}, if 
$s(v+w)\leq s(v)+s(w)$ and $s(tv)=ts(v)$ for $v,w$ in $\Ree(G)$ and $t\geq 0$. 
The set of sublinear functionals is closed under addition, pointwise maximum: $s\vee s'$,
and variable inversion: $\check{s}(v)=s(-v)$.
Aspects of the next proposition can be found in
 \cite[pps. 120--121]{NariciB}, but we have a different proof, whose key inequality
 we shall use in the sequel.

\begin{proposition}\label{prop:sublinear}
If $s$ is sublinear, then $s+\check{s}$ and $s\vee\check{s}$ are seminorms.
Hence, any sublinear functional is continuous on $\fV$.
\end{proposition}

\begin{proof}
We note for $v$ in $\fV$
that $0=s(0)\leq s(v)+s(-v)$, i.e.\ at least one of $s(\pm v)$ is non-negative.
Hence both $s+\check{s}$ and
$s\vee\check{s}$ are sublinear, non-negative, and symmetric, i.e.  seminorms.  

For $v,w$ in $\Ree(G)$ we 
have $s(v)\leq s(v-w)+s(w)$ so $s(v)-s(w)\leq s(v-w)\leq s\vee\check{s}(v-w)$.  Likewise
$s(w)-s(v)\leq s\vee\check{s}(w-v)$, so 
\begin{equation}\label{eq:key-ineq}
|s(v)-s(w)|\leq s\vee\check{s}(v-w).
\end{equation}  Hence
continuity follows from Remark \ref{rem:indtop} (2).
\end{proof}

We set
\[
S(\fV),\quad S_0(\fV),\quad\text{and}\quad S_s(\fV)
\]
to be, respectively, the set of sublinear functionals those which are non-negative, and seminorms.

\begin{remark}\label{rem:sublinreal}
As is well-known, we have
\[
\phantom{m}S(\Ree)=\{s_{a,b}:a\leq b \text{ in }\Ree\}\text{ where }
s_{a,b}(v)=\begin{cases} av &\text{if }v<0 \\ bv &\text{if }0\leq v.\end{cases}
\]
Indeed, if $s\in S(\Ree)$ let $a=s(-1)$ and $b=s(1)$, which, by homogeneity, determines
$s$ on each of $(-\infty,0)$ and $(0,\infty)$, respectively.  Convexity of $s$ gives $a\leq b$.
Notice that
\[
S_0(\Ree)=\{s_{a,b}:a\leq 0\leq b\}\text{ and }S_s(\Ree)=\{s_{-b,b}:b\geq 0\}.
\]
\end{remark}

As in the introduction, we let $G$ be a locally compact abelian group.  We insert the
following result here due to its proof's commonality with Proposition \ref{prop:sublinear}.

\begin{proposition}\label{prop:homogenization}
Let $p:G\to\Ree$ be continuous and subadditive.  Then
\[
p^+(x)=\lim_{n\to\infty}\frac{p(nx)}{n}
\]
defines a continuous, subadditive, and $\En$-homogeneous (i.e.\ $p^+(mx)=mp^+(x)$
for $m$ in $\En$) function on $G$.    Furthermore, $p^+$ is non-negative
if $p$ is bounded below, and symmetric if $p$ is.
\end{proposition}

\begin{proof}
Notice that $p(0)\leq 2p(0)$ so $p(0)\geq 0$, and then for $x$ in $G$
we have $p(0)\leq p(x)+p(-x)$ so $p\vee\check{p}\geq 0$.  If
$\{\frac{p(nx)}{n}\}_{n\in\En}$ is bounded below then Fekete's Lemma
provides that the limit $p^+(x)$ exists, which always happens when $p(x)\geq 0$.
Evidently, $p^+(0)=0$, and if $p$ is bounded below then $p^+\geq 0$.
If $p(x)<0$, then 
\[
-np(-x)\leq -p(-nx)\leq p(0)-p(-nx)\leq
p(nx)\text{ for }n\text{ in }\En
\]
 and we see that $\{\frac{p(nx)}{n}\}_{n\in\En}$ is bounded below by $-p(-x)$.  

All of the desired properties of $p^+$ are straightforward, save for continuity.
For $x,y$ in $G$ and $n$ in $\En$, we use (\ref{eq:key-ineq})
to see that
\[
\left| \frac{p(nx)}{n}-\frac{p(ny)}{n}\right| \leq \frac{1}{n}p\vee\check{p}(n(x-y))\leq 
p\vee\check{p}(x-y).
\]
Thus we see that the family of functions $\{x\mapsto \frac{p(nx)}{n}\}_{n\in\En}$ is 
equicontinuous on $G$.  Hence Arzela-Ascoli theorem shows that 
the limit defining $p^+$ is a uniform limit on compact
subsets of $G$.
\end{proof}

\subsection{Weak* compact, convex sets in the dual}
We give a version of well-known results due to  H\"{o}rmander (\cite{hormander}).
As we wish to emphasize particular aspects,
we have elected to rebuild the proofs from basic results.

We now let $\fV^*$ denote the linear dual and set
\[
\fL(\fV)=\{L\subset \fV^* : L\text{ is non-empty, convex, and weak*-compact}\}
\]
where, of course the weak*-topology is the coarsest topology making
evaluation functionals, $\lam\mapsto \lam(v)$ for $v$ in $\fV$, each continuous.
We shall also distinguish
\[
\fL_0(\fV)=\{L\in\fL(\fV):0\in L\}\text{ and }\fL_s(\fV)=\{L\in\fL(\fV):-L=L\}.
\]
We shall denote the polar of a non-empty subset $B$ of $\fV$ by
\[
B^\circ=\{\lam\in\fV*:\lam(v)\leq 1\text{ for }v\text{ in }B\}
\]
and the prepolar of a non-empty subset $L$ of $\fV^*$ by
\[
L_\circ=\{v\in\fV:\lam(v)\leq 1\text{ for }\lam\text{ in }L\}.
\]

\begin{lemma}\label{lem:sublin_conv}
For $s$ in $S(\fV)$, let
\[
L_s=\{\lam\in\fV^*:\lam\leq s\}.
\]
Then $s\mapsto L_s:S(\fV)\to\fL(\fV)$ is a bijection which takes $S_0(\fV)$ onto $\fL_0(\fV)$
and $S_s(\fV)$, onto $\fL_s(\fV)$.  
\end{lemma}

\begin{proof}
We first show that each $L_s\in\fL(\fV)$.  
 Given $s$ in $S(\fV)$,  and $v$ in $\fV\setminus\{0\}$, Remark \ref{rem:sublinreal} shows that
$s|_{\Ree v}\cong s_{a,b}$, i.e.\ $s(tv)=s_{a,b}(t)$.
We let $\mu:\Ree v\to\Ree$ be given
by $\mu(tv)=ct$ where $a\leq c\leq b$.  Then Hahn-Banach Theorem provides a 
$\lam_0$ in $L_s$ which extends $\mu$, i.e.\ $L_s\not=\varnothing$.  
Next we note for $\lam$ in $L_s$ that $-\lam=\check{\lam}\leq \check{s}$, so $-\check{s}\leq \lam$,
and hence 
\begin{equation}\label{eq:Ls}
L_s=\{\lam\in\fV^*:-\check{s}\leq\lam\leq s\}
\end{equation}
in particular $-\check{s}\leq s$ on $\fV$.
As in proof of Banach-Alaoglu Theorem, the map $\lam\mapsto (\lam(v))_{v\in \fV}$ 
from $L_s$ to the product space
$\prod_{v\in \fV}[-\check{s}(v),s(v)]$ is injective, continuous, and has closed range,
showing that $L_s$ is weak*-compact.  It is evident that $L_s$ is convex; further
$0\in L_s$ if $s\in S_0(\fV)$, and $-L_s=L_{\check{s}}$, so $L_s\in \fL_s(\fV)$
if $s\in S_s(\fV)$.

Let us verify injectivity.  Let us first assume that $s\in S_0(\fV)$.  We consider the ``$s$-ball"
\[
B_s=\{v\in \fV:s(v)\leq 1\}
\]
which is a convex set containing $0$. Now for $\lam$ in $\fV^*$ we have that
$\lam\leq s$ exactly when $\lam(\frac{1}{s(v)}v)\leq 1$ whenever $s(v)>0$, which
is exactly when $\lam\in B_s^\circ$.  Hence $L_s=B_s^\circ$.  The bipolar
theorem  (see, for example, \cite[VI.1.5]{schaeferwolff})
shows that $(L_s)_\circ=B_s$.  Since $s\mapsto B_s$ is injective,
 injectivity of $s\mapsto L_s$ on $S_0(\fV)$ follows.
If we consider $s$ in $S(\fV)$ generally, then for any $\lam_0$ in $L_s$, $s-\lam_0\in S_0(\fV)$
and $L_s=\lam_0+L_{s-\lam_0}$.  If also $s'\in S(\fV)$ then either $L_s\cap L_{s'}=\varnothing$,,
or there is $\lam_0\in L_s\cap L_{s'}$ and $s-\lam_0=s'-\lam_0$ if and only if $L_{s-\lam_0}=
L_{s-\lam_0}$, showing that $s\mapsto L_s$ is injective on $S(\fV)$.

Let us verify surjectivity.  If $L\in \fL(\fV)$ we let
\begin{equation}\label{eq:sublin_fromL}
s(v)=\max\{\lam(v):\lam\in L\}.  
\end{equation}
Evidently, $s\in S(\fV)$; and, further is in $S_0(\fV)$ if $0\in L$, and is in $S_s(\fV)$ if $-L=L$.
If $0\in L$, then for $v$ in $\fV$, $s(v)\leq 1$ exactly when $\lam(v)\leq 1$ for all $\lam$ in $L$,
i.e.\ $B_s=L_\circ$, so, again, the bipolar theorem provides that $L=B_s^\circ$.  Hence, from above,
$L=L_s$.  If $L\in\fL(\fV)$ genrally, then for any $\lam_0$ in $L$, $-\lam_0+L\in\fL_0(\fV)$,
so is $L_{s_0}$ for some $s_0$ in $S_0(\fV)$, hence $L=\lam_0+L_{s_0}=L_{s_0+\lam_0}$.
\end{proof}

\begin{remark}\label{rem:sublin_conv}
\begin{enumerate}
\item As expected, we see that
$\fL(\Ree)=\{[a,b]:a\leq b\text{ in }\Ree\}$
where each $[a,b]=L_{s_{a,b}}$.  Notice for $a\leq 0\leq b$ we have ``balls"
$B_{s_{a,b}}$ of the form $[\frac{1}{a},\frac{1}{b}], \, (-\infty,\frac{1}{b}],\, [\frac{1}{a},\infty),$ or $\Ree$.
\item If $0< a<b$, above, then $B_{s_{a,b}}=(-\infty,\frac{1}{b}]$, giving polar $[0,b]$.
Generally, for $s\in S(\fV)\setminus S_0(\fV)$ there seems to be no reasonable
definition of an associated ball $B_s$ to use a bipolar argument as in the second paragraph
or the proof above.  
\end{enumerate}
\end{remark}

We observe that $\fL(\fV)$ enjoys operation of Minkowski sum,
and convex join:
\[
L+L'=\{\lam+\lam':\lam\in L,\;\lam'\in L'\}\text{ and }
L\vee L'=\conv(L\cup L').
\]
As weak*-topology is a vector space topology, and 
$L\vee L'$ is the range of the map $(t,\lam,,\lam')\mapsto (1-t)\lam+t\lam'$ on
$[0,1]\times L\times L'$, we see that both these sets are compact.
The identity $L_0=\{0\}$ provides that
sum is monoidal, and that the join is monoidal when restricted to $\fL_0(\fV)$, hence on
$\fL_s(\fV)$.
We consider operations on $\fL(\fV)$ as witnessed by operations on $S(\fV)$.

\begin{lemma}\label{lem:opsonLV}
Let $s,s'\in S(\fV)$, and $T:\fV'\to\fV$ be a linear operator from another real vector space.
Then we have
\[
L_{s\vee s'}=L_s\vee L_{s'},\quad L_s+L_{s'}=L_{s+s'},\quad\text{and}\quad
T^*(L_s)=L_{s\circ T}
\]
where $T^*:\fV^*\to{\fV'}^*$ is the adjoint.
\end{lemma}

\begin{proof}
Lemma \ref{lem:sublin_conv} provides $s''$ in $S(\fV)$ for which $L_s\vee L_{s'}=L_{s''}$.
Then for $\lam''$ in $\fV^*$ we have
\begin{align*}
\lam''\leq s'' \quad
&\Leftrightarrow\quad \lam''=(1-t)\lam+t\lam'\text{ where }
(t,\lam,\lam')\in [0,1]\times L_s\times L_{s'} \\
&\Rightarrow\quad \lam''\leq s\vee s'
\end{align*}
showing that $L_{s''}\leq L_{s\vee s'}$, i.e.\ $s''\leq s\vee s'$.  But $L_s\cup L_{s'}\subseteq
L_{s''}$, so $s,s'\leq s''$, hence $s\vee s'\leq s''$.

The proof for the sum is similar to, if mildly simpler than, that of the join, and we omit it.

For the linear operation, first suppose that $s\in S_0(\fV)$. Then as in the 
the proof of Lemma \ref{lem:sublin_conv}
we have $L_s=B_s^\circ$.  Thus
\[
T^*(L_s)=T^*(B_s^\circ)=T(B_s)^\circ=B_{s\circ T}^\circ=L_{s\circ T}.
\]
If, generally, $s\in S(\fV)$, fix $\lam_0$ in $L_s$ so $L_s=\lam_0+L_{s-\lam_0}$
where $s-\lam_0\in S_0(\fV)$.  Then
\[
T^*(L_s)=T^*(\lam_0+L_{s-\lam_0})=\lam_0\circ T+L_{(s-\lam_0)\circ T}=L_{s\circ T}.
\]
\end{proof}

\begin{remark}
\begin{itemize}
\item Notice that $\fV^*\subseteq S(\fV)$, and for
$\lam_0$ in $\fV^*$ that $L_{\lam_0}=\{\lam_0\}$.
Hence the sum and linear action formulas above show the effect of weak*-continuous affine
transformations on $\fL(\fV)$.
\item If $L_s\cap L_{s'}\not=\varnothing$, then 
$L_s\cap L_{s'}=L_{s\wedge s'}$ where $s\wedge s'(v)=\max\{\lam(v):\lam\in L_s\cap L_{s'}\}$.
This is immediate from (\ref{eq:sublin_fromL}).
\item The semilattice $(\fL(\fV),\vee)$ is Dedekind complete: if there is $L$ in $\fL(\fV)$
so each $L_i\subseteq L$, then $\bigvee_{i\in I}L_i=\wbar{\conv}^{w*}\bigcup_{i\in I}L_i\in\fL(\fV)$.
The semilattice $(\fL(\fV)\cup\{\varnothing\},\cap)$ is complete, as are sublattices
$\fL_0(\fV)$ and $\fL_s(\fV)$.
 The lattice $(\fL(\fV)\cup\{\varnothing\},\vee,\cap)$ is not distributive, even for $\fV=\Ree$,
 but sublattices $\fL_0(\fV)$ and $\fL_s(\fV)$ are distributive.  Analogous
results carry to $S(\fV)$.
\end{itemize}
\end{remark}

\subsection{The realization of a locally compact abelian group}
To begin, let the {\em (topological) torsion} subgroup, also known as the group
of compact elements, be given by
\[
\tor(G)=\left\{x\in G:\wbar{\{kx\}_{k\in\Zee}}\text{ is compact}\right\}
\]
which is a closed subgroup by \cite[(9.10)]{hewittrI}.
[Though in many cases $\tor(G)$ is a factor of $G$, this is not always
the case; see Section \ref{ssec:nonsplit}.]
Then $G/\tor(G)$ admits no non-trivial compact subgroups, so it dual
$\what{G/\tor(G)}$ admits no open subgroups, and thus is connected.
Hence, by \cite[2.4.1]{rudin}, for example, we find that $\what{G/\tor(G)}$
is isomorphic to
$\Ree^d\times T$, where $T$ is a compact connected group.
Thus Pontryagin duality provides that
\[
G/\tor(G)\cong \Ree^d\oplus\Gam
\]
where $\Gam=\what{T}$ is a torsion-free discrete group.  We let
$q:G\to G/\tor(G)$ denote the quotient map so
\[
q(x)=q_0(x)\oplus q_1(x)\text{ where }q_0(x)\in\Ree^d\text{ and }q_1(x)\in\Gam.
\]
We define the {\em realization} of $G$ by 
\[
\Ree(G)=\Ree^d\oplus(\Ree\otimes\Gam)
\]
with tensoring over $\Zee$.  This is a real vector space.  The {\em rank} of $G$
is given by
\[
\rank G=\dim\Ree(G).
\]
We let $\iota:G\to \Ree(G)$ be given by 
\[
\iota(x)=q_0(x)\oplus(1\otimes q_1(x)).
\]
We consider on $\Ree(G)$ the finest locally convex topology,  discussed in 
Section \ref{ssec:flctop}  Thus
\begin{equation}\label{eq:iota_cts}
\ker\iota=\tor(G)\text{ and }\iota:G\to\Ree(G)\text{ is continuous}
\end{equation}
where the latter holds as $q_0$ is continuous.

We denote families of sublinear functionals
\[
S(G)=S(\Ree(G)),\quad S_0(G)=S_0(\Ree(G)),\quad\text{and}\quad S_s(G)=S_s(\Ree(G)).
\]
In the case $G=\Ree^d$ we note that $\Ree(G)=\Ree^d$, so there will be no ambiguity in the notation
$S(\Ree^d)$.

The following may be regarded as the {\em universal property} for $\Ree(G)$.
Recall that $\En$-homogeneity is defined in Proposition \ref{prop:homogenization}.

\begin{proposition}\label{prop:extension}
\begin{enumerate}
\item Let $p:G\to\Ree$ be continuous, subadditive, and $\En$-homogeneous.
Then there is a unique $s$ in $S(G)$ for which $s\circ\iota=p$, and each
$s$ in $S(G)$ is realized this way.

Furthermore, if $p':G\to\Ree$ is also continuous, subadditive, and $\En$-homogeneous,
with $s'$ in $S(G)$, as above, then $p'\leq p$ if and only if $s'\leq s$.  Hence
$p\geq 0$ implies $s\in S_0(G)$ and $p=\check{p}$ (where $\check{p}(x)=p(-x)$) implies $s\in S_s(G)$.
\item If $\eta:G\to\fW$ is a continuous additive map into a locally convex vector
space, then there is a unique linear $\eta_\Ree:\Ree(G)\to\fW$ for which
$\eta_\Ree\circ\iota=\eta$.
\end{enumerate}
\end{proposition}

\begin{proof}
(1) We have for $y$ in $\tor(G)$ that $p(y)=\lim_{n\to\infty}\frac{p(ny)}{n}=0$, so for
$x$ in $G$ we have
\[
p(x+y)\leq p(x)+p(y)=p(x)\leq p(x+y)+p(-y)=p(x+y)
\]
i.e.\ $p(x+y)=p(x)$.  Hence $\bar{p}:G/\tor(G)\cong\Ree^d\oplus\Gam\to\Ree$
given by $\bar{p}\circ q=p$ is well defined and continuous, and easily checked
to be $\En$-homogeneous and subadditive.

We consider the ``rationalization" of $G$, $\Que(G)=
\Ree^d\oplus(\Que\otimes\Gam)$ which is a $\Que$-vectors space and contains $\iota(G)$.
We define on $\Que(G)$
\begin{equation}\label{eq:rational}
\til{p}\left(x,\sum_{j=1}^k \frac{m_j}{n_j}\otimes y_j\right)=
\frac{1}{n}\bar{p}\left(nx,\sum_{j=1}^k \frac{nm_j}{n_j} y_j\right)
\end{equation}
where $x\in\Ree^d$, for each $j$,  $y_j\in\Gam$, $m_j\in\Zee$, $n_j\in\En$, and
$n\in\En$ is divisible by $n_1\dots n_k$.  It is easily checked that
$\til{p}$ is well-defined, satisfies $\til{p}(\frac{m}{n}v)=\frac{m}{n}\til{p}(v)$
for $m,n$ in $\En$ and $v$ in $\Que(G)$, and is subadditive, 
i.e.\ $\til{p}$ is rationally sublinear.  Finally, we have that $\til{p}\circ\iota=\bar{p}\circ q=p$,

To obtain $s$, let $v=\left(x,\sum_{j=1}^kr_j\otimes y_j\right)$ in $\Ree(G)$, 
where $x\in\Ree^d$, each $r_j\in\Ree$ and 
each $x_j\in \Gam$.  We may further assume that $\{1\otimes y_1,\dots,1\otimes y_k\}$ is linearly 
independent in $\Ree\otimes\Gam$.  For each $j$ let $(q_{jl})_{l=1}^\infty$ be a sequence of 
rationals tending to $r_j$, and define
\[
s(v)=\lim_{l\to\infty}\til{p}(v_l)\text{ where each }v_l=\left(x,\sum_{j=1}^kq_{lj}\otimes y_j\right)
\]
Indeed it follows from the evident analogue of (\ref{eq:key-ineq})  that
\begin{align*}
|\til{p}(v_l)-\til{p}(v_{l'})|\leq \til{p}\vee\check{\til{p}}(v_l-v_{l'}) 
\leq \sum_{j=1}^k |q_{lj}-q_{l'j}|\til{p}\vee\check{\til{p}}(0,1\otimes y_j)
\end{align*}
showing that $(\til{p}(v_l))_{l=1}^\infty$ is Cauchy, hence the limit $s(v)$ exists.  Likewise,
this limit is independent of choices of rational sequences $(q_{jl})_{l=1}^\infty$ tending to respective 
$r_j$.   The sublinearity of $s$ is immediate from the rational sublinearity of $\til{p}$.
Moroever, $s$ extends $\til{p}$, so $s\circ\iota=\til{p}\circ\iota=p$.

Notice that if $s''\in S(G)$ with $s''\circ\iota=p$, then $s''=\til{p}$ on the dense set
$\Que(G)$ of $\Ree(G)$.  By continuity, this forces $s''=s$, and hence this $s$ is unique.

Finally, $p'\leq p$ exactly when $\til{p}'\leq\til{p}$, which, again by continuity, is exactly when
$s'\leq s$.  Hence non-negativity of $p$ begets that of $s$.  That symmetry
of $p$ begets that of $s$ is similar.

(2) As above, we notice that $\tor(G)\subseteq \ker\eta$, and we may evidently
define additive continuous $\bar{\eta}:G/\tor(G)\to\fW$ with $\bar{\eta}\circ q=\eta$,
and $\bar{\eta}(mz)=m\bar{\eta}(z)$ for
$m$ in $\Zee$ and $z$ in $G/\tor(G)$.  Hence we can define
$\til{\eta}:\Que(G)\to\fV$ by
\[
\til{\eta}\left(x,\sum_{j=1}^k \frac{m_j}{n_j}\otimes y_j\right)
=\frac{1}{n}\bar{\eta}\left(x,\sum_{j=1}^k \frac{nm_j}{n_j}y_j\right)
\]
where data is as in (\ref{eq:rational}), above, so $\til{\eta}$ is $\Que$-linear.  
The proof of (1) is easily applied to see that $\til{\eta}$ extends uniquely to a linear 
$\eta_\Ree:\Ree(G)\to\fW$, with $\eta_\Ree\circ\iota=\eta$.
\end{proof}

\begin{remark}\label{rem:homo}
Let $G'$ be another locally compact abelian group with $\iota':G'\to\Ree(G')$ the cannonical
map into its realization, and let $\alp:G'\to G$ be a continuous homomorphism.
Then $\iota\circ\alp:G'\to\Ree(G)$ is a continuous and additive, and we let
$\alp_\Ree=(\iota\circ\alp)_\Ree:\Ree(G')\to\Ree(G)$ so $\alp_\Ree\circ\iota'=\iota\circ\alp$.
\end{remark}

We let $\Ree(G)^*$ denote the (continuous) linear dual of $\Ree(G)$.  
The set $\{\lam\circ\iota:\lam\in \Ree(G)^*\}$ is the set of {\em real characters}, defined
in \cite[(24.33)]{hewittrI}. We let
$uc$ denote the topology of uniform convergence on compact sets.

\begin{lemma}\label{lem:wstar_uonc}
Let $(\lam_i)_i$ be a net in $\Ree(G)^*$.  Then weak*-$\lim\lam_i=\lam$ in $\Ree(G)^*$
if and only if $uc$-$\lim\lam_i\circ\iota=\lam\circ\iota$ on $G$.
\end{lemma}

\begin{proof}
Notice that if  weak*-$\lim\lam_i=\lam$, then norm-$\lim\lam_i|_F=\lam|_F$ for any finite dimensional
subspace $F$.  Any compact
subset $K$ of $G$ has that $\iota(K)\subseteq L\times E$ where $L$ is compact in
$F=\iota(\Ree^d)$, where $\Ree^d$ is
the connected component of the identity of $G/\tor(G)$, and $E$ is finite.  

Conversely, if $uc$-$\lim\lam_i\circ\iota=\lam\circ\iota$, then convergence is pointwise, and 
$\iota(G)$ spans $\Ree(G)$.
\end{proof}

\section{Spectra of Beurling algebras}\label{sec:spectrum}

\subsection{The complex dual group}\label{ssec:complexification}
We let
\begin{align*}
\what{G}_\Cee&=\left\{\sig:G\to\Cee\setminus\{0\}\;\left| \;\begin{matrix}
\sig\text{ is continuous with} \\ \sig(x+y)=\sig(x)\sig(y)\text{ for }x,y\text{ in }G\end{matrix}
\right\}\right.  \\
\what{G}\phantom{f}&=\{\sig\in\what{G}_\Cee\;|\;\sig(G)\subseteq\Tee\} \\
\what{G}^+&=\{\sig\in\what{G}_\Cee\;|\;\sig(G)\subseteq (0,\infty)\}
\end{align*}
Notice that $\what{G}_\Cee$ is a group under pointwise multiplication, and
elements of $\what{G}_\Cee$ enjoy a polar decomposition, $\sig=\frac{\sig}{|\sig|}|\sig|$,
showing that $\what{G}_\Cee=\what{G}\what{G}^+$.  

Under topology of uniform convergence on compact subsets, $\what{G}_\Cee$
is may be easily checked to be a topological group.  Indeed if
$uc$-$\lim_i\sig_i=\sig$ in $\what{G}_\Cee^+$ then on any compact set $K$
the characters $\sig_i$ are eventually bounded; hence the continuity of products and inverses
follows.  It is straightforward to see that $\what{G}^+$ and $\what{G}$ are closed subgroups
and polar decomposition provides an isomorphism of topological groups:
\[
\what{G}_\Cee\cong \what{G}\times\what{G}^+.
\]

\begin{proposition}\label{prop:poschar}
The map $E:\Ree(G)^*\to\what{G}^+$ given by 
\[
E(\lam)=\exp\circ\lam\circ\iota
\]
is a bijection and a homeomorphism for the weak*-topology on $\Ree(G)^*$ and
the topology of uniform convergence on compact sets on $\what{G}^+$.  Hence
$\what{G}^+$ is complete.
\end{proposition}

\begin{proof}
Proposition \ref{prop:extension} shows that $\lam\mapsto\lam\circ\iota$ is injective,
so $E$ is injective.  Given $\sig\in\what{G}^+$, the same proposition shows that
$\lam=(\log\circ\sig)_\Ree\in\Ree(G)^*$ with $E(\lam)=\sig$.  The
homeomorphism property is immediate from Lemma \ref{lem:wstar_uonc}.
Note that a dual space is complete (with respect to additive group uniformity).
\end{proof}

As a consequence we see that
\[
 \what{G}^+= \what{G/\tor(G)}^+\circ q
\]
and this set separates points of $G/\tor(G)$ since linear functionals separate points
of $\Ree(G)$.  

\begin{remark}
We get other natural functorial properties for a closed subgroup $H$ of
$G$:
\begin{equation}\label{eq:functorial}
\what{G/H}^+\hookrightarrow \what{G}^+\text{ and }
\what{G}^+|_H=\what{H}^+.
\end{equation}
The first is obvious by taking adjoint of the quotient map. The second follows as
$\tor(H)=\tor(G)\cap H$, so $\iota_G|_H$ extends to an isomorphism of $\Ree(H)$ onto its range
$\spn \iota_G(H)$ of $\Ree(G)$.  Hence since all linear functional extend,
$\Ree(G)^*|_{\Ree_G(H)}=\spn \iota_G(H)^*$.  Thus $\what{H}^+
=\{\exp\circ\lam\circ\iota_G|_H:\lam\in\Ree(G)^*\}=\what{G}^+|_H$.  Compare this with
\cite[(24.36)]{hewittrI}.
\end{remark}

Let
\[
\fK=\{K\subseteq G\;|\; K\text{ is compact with }K^\circ\not=\varnothing\}
\]
where $K^\circ$ is the interior of $K$.  For each $K$ we identify $\bl^1(K)\cong 1_K\bl^1(G)$, where
$\bl^1(G)$ is the usual group algebra consisting of (equivalence classes of)
complex Haar integrable functions on $G$, and let
\begin{equation}\label{eq:lonec}
\bl^1_c(G)=\bigcup_{K\in\fK}\bl^1(K).
\end{equation}
We assign to $\bl^1_c(G)$ the inductive topology, which is the finest locally
convex topology allowing each inclusion $\bl^1(K)\hookrightarrow \bl^1_c(G)$ to be continuous.
This topology has neighbourhood basis about $0$ consisting of convex, balanced, absorbing
sets $V$ for which $V\cap\bl^1(K)$ is open in $\bl^1(K)$ for each $K$ in $\fK$.  See 
\cite[I.6]{schaeferwolff}.

\begin{remark} \label{rem:contoncomp}
If $\fW$ is a complex locally convex space, in particular a Banach space, and 
$T:\bl^1_c(G)\to\fW$ is linear with each $T|_{\bl^1(K)}$ continuous, then
$T$ is continuous.  Hence, this topology admits the following facts, where $K,K'\in\fK$.
\begin{enumerate}
\item The continuous dual of $\bl^1_c(G)$ is the space
$\bl^\infty_{loc}(G)$ of (equivalence classes) of Borel measurable functions
which are essentially bounded on each $K$, via the usual
dual pairing: 
\[
\langle f,\vphi\rangle=\int_G f\vphi
\]
for $f$ in $\bl^1(G)$ and $\vphi$ in $\bl^\infty_{loc}(G)$.  
\item Given $f$ in $\bl^1(K)$, the map $g\mapsto f\ast g:\bl^1(K')\to\bl^1(K+K')$ is continuous
and hence extends to a continuous operator on $\bl^1_c(G)$.  Likewise, given
$x$ in $G$, the map $g\mapsto x\ast g:\bl^1(K')\to\bl^1(x+K')$, where
$x\ast g(y)=g(y-x)$ for a.e.\ $y$, is continuous; and $y\mapsto y\ast g:U\to L^1(\wbar{U}+K')$
is continuous, where $U$ is any relatively compact neighbourhood of $x$. 
Hence $\bl^1_c(G)$ is an (semi-)topological
algebra under convolution which admits continuous translation.
\end{enumerate}
\end{remark}

Let us quickly review the main relationship between the two objects introduced.  We let
\[
\Sig(\bl^1_c(G))=\left\{\sig\in\bl^\infty_{loc}(G)\setminus\{0\}\;\left|\; 
\begin{matrix} \langle f\ast g,\sig\rangle=\langle f,\sig\rangle\langle g,\sig\rangle \\
\text{ for }f,g\text{ in }\bl^1_c(G)\end{matrix}\right\}\right. .
\]
The following can be shown by making evident adaptations of standard proofs (see, for example,
Theorems 2.8.2 and 2.8.5 of \cite{kaniuth}), with appropriate continuity results provided by
Remark \ref{rem:contoncomp}.

\begin{proposition}\label{prop:lonecspec}
We have that $\Sig(\bl^1_c(G))=\what{G}_\Cee$.  The weak* topology on $\what{G}_\Cee$
coincides with the topology $uc$ of uniform convergence on compact sets. 
\end{proposition}

\begin{corollary}\label{cor:homeo}
With weak*-topology, $\what{G}^+$ a closed subgroup of $\what{G}_\Cee$.
Hence both of $\what{G}_\Cee$ and $\what{G}^+$ are 
locally compact if and only if $\rank G<\infty$.
\end{corollary}

\begin{proof}
We use Proposition \ref{prop:poschar}.
Since $\what{G}^+$ is $uc$-closed, it is weak*-closed.
By, for example \cite[I.3.6]{schaeferwolff},
 $\Ree(G)^*$ is locally compact if and only if it is finite dimensional, hence
exactly when $\rank G=\dim\Ree(G)<\infty$.  The characterization of local compactness
of $\what{G}^+$, or $\what{G}_\Cee$, follows.
\end{proof}

\subsection{The spectra of Beurling algebras}
We let $\ome\in W(G)$.  Then for $x$ in $G$ define
its {\em homogenization}
\[
\ome^+(x)=\lim_{n\to\infty}\ome(nx)^{1/n}
\]
where it is is well-known that this limit exists, and  infemum of $\ome(nx)^{1/n}$;
see \cite[A.1.2.6]{dales}, for example.  Being a squential limit of continuous functions, 
$\ome^+$ is Borel measurable, and easily seen to also be a weight.
However, $\ome^+$  is continuous.

\begin{proposition}\label{prop:se-form}
If $\ome\in W(G)$ then there is a unique sublinear functional
$s_\ome$ on $\Ree(G)$ for which $\ome^+=\exp\circ s_\ome\circ\iota$.
\end{proposition}

\begin{proof}
Let $p=\log \ome$, and $p^+$ its homogenization which is continuous by
 Proposition \ref{prop:homogenization}.
Then Proposition \ref{prop:extension} provides a unique sublinear functional
$s_\ome$ with $s_\ome\circ\iota=p^+$.  
\end{proof}

\begin{remark} 
We note that if we assume only that $\ome$ is Borel measurable, then as
noted in Remark \ref{rem:borel} it is equivalent to continuous weight, so $\ome^+$ is
continuous.
\end{remark}

We define the sets of {\em homogeneous weights} by
\[
W^+(G)=\{\ome\in W(G):\ome(nx)=\ome(x)^n\text{ for }x\text{ in }G\text{ and }n\text{ in }\En\}
\]
and analogously define
\[
W^+_0(G)=W_0(G)\cap W^+(G)\text{ and }W^+_s(G)=W_s(G)\cap W^+(G).
\]
Notice that $\what{G}^+\subseteq W^+(G)$ with $\what{G}^+\cap W^+_0(G)=\{1\}$.
The next proposition is an extension of Proposition \ref{prop:poschar} which admits the same proof.

\begin{proposition}\label{prop:homweight-char}
The map $E:S(G)\to W(G)$ given by
\[
E(s)=\exp\circ s\circ\iota
\]
is injective, with 
\[
E(S(G))=W^+(G), \; E(S_0(G))=W_0^+(G),\text{ and }
E(S_s(G))=W_s^+(G).
\]
\end{proposition}

We now consider the Gelfand spectra of Beurling algebras $\bl^1(G,\ome)$.
If $\ome\in W(G)$, let 
\[
\what{G}_\ome=\{\sig\in\what{G}_\Cee:|\sig|\leq\ome\}.
\]
Polar decomposition shows that $\what{G}_\ome=\what{G}\times\what{G}_\ome^+$
where $\what{G}_\ome^+=\what{G}_\ome\cap\what{G}^+$.
Part (3), below, is the culminating result and entire motivation for this note.

\begin{theorem}\label{theo:BAspec}
Let $\ome\in W(G)$.
\begin{enumerate}
\item The Gelfand spectrum $\Sig(\bl^1(G,\ome))$ is isomorphic
to $\what{G}_\ome$.
\item There is $L_\ome\in\fL(G)$ for which $E(L_\ome)=\what{G}_\ome^+$,
so $\what{G}_\ome^+$ is weak*-compact.  The decomposition 
$\what{G}_\ome=\what{G}\times\what{G}_\ome^+$ is a homeomorphism.
\item The map $\ome\mapsto L_\ome:W(G)\to \fL(G)$ is surjective, and bijective
when restricted to $W^+(G)$.  Furthermore, $W_0(G)$ is mapped to $\fL_0(G)$,
while $W_s(G)$ is mapped to $\fL_s(G)$.
\end{enumerate}
Hence all possible spectra of Beurling algebras on $G$ are parameterized by $\fL(G)$, i.e.\
it is the family of all subsets of $\what{G}_\Cee$ given by 
\[
\{\what{G}E(L):L\in\fL(G)\}.
\]
\end{theorem}

\begin{proof}
(1) This is well-known, as we see, for example, by \cite[2.8.2]{kaniuth}.  Let us note, however,
that it may be seen as 
a combination of three facts.  First, $\bl^1_c(G)$ is dense in $\bl^1(G,\ome)$,
so Proposition \ref{prop:lonecspec} provides that $\Sig(\bl^1(G,\ome))$ may be identified as a subset
of $\what{G}_\Cee$.  Second, the dual of $\bl^1(G,\ome)$ consists of elements $\vphi$ of
$\bl^\infty_{loc}(G)$ for which $\frac{\vphi}{\ome}$ is essentially bounded.
Third, characters on a Banach algebra are automatically contractive.

(2) Consider $\sig\in\what{G}^+$.  Then if $\sig\leq\ome$, we find $\sig=\sig^+\leq\ome^+$.
Propositions \ref{prop:se-form} and \ref{prop:homweight-char}
provide $s_\ome$ in $S(G)$ for which $\ome^+=E(s_\ome)$.  Letting $L_\ome=L_{s_\ome}$
(with notation as in Proposition \ref{prop:se-form} and Lemma \ref{lem:sublin_conv}),
Proposition \ref{prop:poschar}, and the inequalities of Proposition \ref{prop:extension}
provide that $E(L_\ome)=\what{G}^+_\ome$.  Corollary \ref{cor:homeo} then provides that
$\what{G}^+_\ome$ is weak*-compact, and provides the homeomorphism.

(3) For $s\in S(G)$, it is obvious that $\ome=\exp\circ s\circ\iota\in W(G)$,
with $s_\ome=\ome$.  Hence Lemma \ref{lem:sublin_conv} provides the surjectivity.
Furthermore, we see that
$s_{\ome^+}=s_\ome$, so surjectivity carries to subdomain $W^+(G)$.
The uniqueness of Proposition 2.9 provides 
for distinct $\ome,\ome'$ in $W^+(G)$ that $s_\ome\not=s_{\ome'}$, hence
giving injectivity on $W^+(G)$.  If $\ome\in W_0(G)$, then Propositions \ref{prop:homogenization}
and \ref{prop:se-form} provide that $s_\ome\in S_0(G)$, hence Lemma \ref{lem:sublin_conv}
shows that $L_\ome\in\fL_0(G)$.  The preservation of symmetry is similar.
\end{proof}

\begin{remark}
Equation (\ref{eq:Ls}) shows for $\ome$ in $W(G)$ that
\[
\what{G}^+_\ome=E(L_\ome)=
\left\{\sig\in \what{G}^+\;\left|\; \frac{1}{\check{\ome}^+}\leq \sig\leq \ome^+
\right\}\right. .
\]
This is worked out explicitly for $G=\Ree,\Zee$ in many sources; see, for example, Propositions
2.8.6 and 2.8.8 in \cite{kaniuth}, and compare with Sections \ref{ssec:reals} and
\ref{ssec:rational}, below.
\end{remark}

We note for $\ome,\ome'$ in $W(G)$ we have that $\ome\ome',\ome\vee\ome'\in W(G)$.
Furthermore, these operations preserve $W_0(G)$ and $W_s(G)$, and homogenization:
$(\ome\ome')^+=\ome^+\ome'^+,(\ome\vee\ome')^+=\ome^+\vee\ome'^+$.  Hence we
obtain to following observation.

\begin{remark}\label{rem:retract}
The homogenization map $\ome\mapsto\ome^+:W(G)\to W^+(G)$ is a semigroup retract
for both multiplication and the join operation.
\end{remark}

\begin{corollary}
For $\ome,\ome'$ in $W(G)$ we have
\[
L_{\ome\ome'}=L_\ome+L_{\ome'}\text{ and }L_{\ome\vee\ome'}=L_\ome\vee L_{\ome'}.
\]
Hence $\what{G}_{\ome\ome'}=\what{G}_\ome\what{G}_{\ome'}$ in $\what{G}_\Cee$.
\end{corollary}

\begin{proof}
It follows Remark \ref{rem:retract} that $s_{\ome\ome'}=s_\ome+s_{\ome'}$ and
$s_{\ome\vee\ome'}=s_\ome\vee s_{\ome'}$.  We then appeal to Lemma \ref{lem:opsonLV}.
Then 
\[
\what{G}_{\ome\ome'}=\what{G}E(L_{\ome\ome'})=\what{G}E(L_\ome+L_{\ome'})
=\what{G}E(L_\ome)E(L_{\ome'})=\what{G}_\ome\what{G}_{\ome'}
\]
as $\what{G}=\what{G}\what{G}$.
\end{proof}

We consider the Gelfand-Raikov-Shilov weights
\[
W_{GRS}(G)=\{\ome\in W(G):\ome^+=1\}
\]
which is evidently also stable for multiplication and the join operation.  Hence
for $\ome\in W(G)$, $\what{G}_\ome=\what{G}$ if and only if $\ome\in W_{GRS}(G)$.

\begin{corollary}\label{cor:GRS_minimal_weights}
\begin{enumerate}
\item For $\ome\in W(G)$, $\ome^+\in \what{G}^+$ if and only if $\ome\in \what{G}^+W_{GRS}(G)$.
In this case $\bl^1(G,\ome)$ is isomorphic to $\bl^1(G,\ome')$ where $\ome'\in W_{GRS}(G)$.
\item For $\ome\in W(G)$, there is $\ome'\in W_0(G)$ for which 
$\bl^1(G,\ome)$ is isomorphic to $\bl^1(G,\ome')$
\item The minimal elements of the semilattice $(W(G),\vee)$ are exactly the elements of $\what{G}^+$.
Hence the corresponding Beurling algebra are all isomorhic to $\bl^1(G)$ (i.e.\ with weight $1$).
\end{enumerate}
\end{corollary}

\begin{proof}
(1)  If $\ome^+=\sig\in\what{G}^+$, then $\ome'=\ome\sig^{-1}\in W_{GRS}(G)$.  
Hence $f\mapsto f\sig$ is an isomorphism from $\bl^1(G,\ome)$ onto $\bl^1(G,\ome')$.

(2) Here we let $\sig\in \what{G}_\ome^+$ and let $\ome'=\ome\sig^{-1}$.  Note that
$L_{\ome\sig^{-1}}=L_\ome-\lam\in\fL_0(G)$ where $\sig=E(\lam)$.  We use the algebra isomorphism
from (1), above.

(3) Minimal elements in $(\fL(G),\vee)$ are singletons, so minimal elements of 
$W^+(G)$ are elements of $\what{G}^+$.  Since the limit defining $\ome^+$ is an infemum,
$\ome\geq \ome^+$.  Hence by Remark \ref{rem:retract}, $\what{G}^+$ is the set of
minimal elements of $W(G)$.
\end{proof}

\begin{corollary}
Let $G'$ be another locally compact abelian group and $\alp:G'\to G$ be a continuous
homomorphism.  Then for $\ome$ in $W(G)$, $\ome\circ\alp\in W(G')$.
Moroeover, $L_{\ome\circ\alp}=\alp_\Ree^*(L_\ome)$ where
$\alp_\Ree$ is defined in Remark \ref{rem:homo}.
\end{corollary}

\begin{proof}
The first statement is obvious.
Since $L_\ome=L_{s_\ome}$, as in part (2) of the proof of Theorem \ref{theo:BAspec},
Lemma \ref{lem:opsonLV} tells us that $\alp_\Ree^*(L_\ome)=L_{s_\ome\circ\alp_\Ree}$.
Hence, it suffices to show that
$s_{\ome\circ\alp}\circ\iota'=s_\ome\circ\alp_\Ree\circ\iota'$, where $\iota':G'\to\Ree(G')$ is
the realization map.  We have for $x$ in $G'$ that
\[
s_\ome\circ\alp_\Ree\circ\iota'(x)=(\log\ome^+)\circ\alp(x)
=\log(\ome\circ\alp)^+(x)=s_{\ome\circ\alp}\circ\iota'(x).
\]
\end{proof}

\begin{remark}
\begin{itemize}
\item We note for any of the subsets of $W(G')$ defined in this article, e.g.\
$W_0(G'),W_s(G'),\what{G'}^+,W_{GRS}(G'),W^+(G')$, we have that
$\ome\mapsto\ome\circ\alp$ carries this set into its evident counterpart, and respects
multiplication and the join.
\item If $\alp$ is an open homomorphism with $H=\alp(G')$, then $\alp$ can be factored as
$\alp=\bar{\alp}\circ q'$, where $q':G'\to G/\ker\alp$ is the quotient map
and  $\bar{\alp}$ is an automorphism on $H$.
It is then possible to devise a surjective homomorphism $\alp_1:\bl^1(G',\ome\circ\alp)\to
\bl^1(H,\ome|_H)$ by
\[
\phantom{mmm}\alp_1(f)=\int_H \bar{f}\circ\bar{\alp}^{-1}\text{ where }\bar{f}(x+\ker\alp)=\int_{\ker\alp}f(x+y)\,dy
\]
where $f\mapsto \bar{f}$ can be made a quotient map   
when Haar measures are appropriately normalized, and $\alp_1$ is an isometry only if
$\alp$ is unimodular.  
Otherwise we expect only to devise a homormorphism into a certain measure Beurling
algebra; compare, for example, with \cite{gharamani_zadeh}.
\end{itemize}
\end{remark}

\section{Examples}\label{sec:examples}

\subsection{The real line}\label{ssec:reals}
Here $\Ree(\Ree)\cong\Ree$.  Proposition \ref{prop:poschar} gives 
$\what{\Ree}^+=\{x\mapsto e^{ax}:a\in\Ree\}\cong\Ree$, $S(\Ree)$ is described in Remark
\ref{rem:sublinreal}, and $\fL(\Ree)$ is described in Remark \ref{rem:sublin_conv}.

\subsection{Rational groups}\label{ssec:rational}
Let $\Gam$ be any  non-zero subgroup of $\Que$.  Then a computation akin to
(\ref{eq:rational}) shows the
well-known fact that $\Que(\Gam)=\Que\otimes \Gam\cong \Que$, and hence $\Ree(G)
\cong\Ree$.  Hence $\what{\Gam}^+$, $S(\Gam)$ and $\fL(\Gam)$ are described exactly as
for $\Ree$, above.

\begin{remark}\label{rem:cotorsion}
\begin{enumerate}
\item Let $\Gam$ be discrete and $\Lambda$ a subgroup with $\Gam/\Lambda$ a torsion
group.  Then $\rank\Lambda=\rank\Gam$.  Indeed, for $x$ in $\Gam$ with $nx\in\Lambda$
we have $1\otimes x=\frac{1}{n}\otimes nx\in\Ree(\Lambda)$, from which it follows that
$\Ree(\Gam)=\Ree(\Lambda)$.
\item If $T$ is compact, and $D$ is a totally disconnected closed subgroup,
then $\rank\what{T/D}=\rank\what{T}$.  Indeed, $\what{T/D}\cong D^\perp\subseteq\what{T}$, and
$\what{T}/D^\perp\cong\what{D}$ is torsion.
\item If $T$ is compact, with connected component $T_e$, then $\rank\what{T}=\rank\what{T_e}$.
Indeed, as above, $\what{T}/T_e^\perp\cong \what{T_e}$ where $T_e^\perp\cong \what{T/T_e}$
is a torsion group.
\end{enumerate}
\end{remark}

\begin{remark}\label{rem:lie}
Any non-zero subgroup $\Gam$ of $\Que$ is of rank 1, and in agreement with Corollary \ref{cor:homeo},
we have that $\what{\Gam}_\Cee\cong \what{\Gam}\times\Ree$ is locally compact.
For any non-cyclic such $\Gam$, $\what{\Gam}$ is a non-Lie solenoid.

In \cite[Theo.\ 5]{CartwrightMc}, it is suggested that $\what{\Gam}$ must be Lie for 
$\what{\Gam}_\Cee$ to be locally compact. However there is an error in their proof regarding
the structure of finite rank abelian groups, as shown by solenoids.
We propose the following correction to \cite[Theo.\ 5]{CartwrightMc}.

{\it For a (not necessarily abelian) compact group $K$ the following are equivalent.
\begin{itemize}
\item[(a)] $K_\Cee$ is locally compact.
\item[(c)] For the connected component $K_e$ and its centre $Z$, $K_e/Z$ is Lie and
$\rank\what{Z}<\infty$.
\end{itemize}}

The omitted condition (b) has to do with a certain space of derivations, which 
we address, below.  Using the reductions of
their paper, we see that $K_\Cee$ is locally compact exactly when
$(K_e)_\Cee$ is locally compact.  A noted in  \cite[Theo.\ 4]{CartwrightMc}, if we
realize $K_e\cong (T\times S)/C$, where $T$ is connected abelian, $S$ is a product of
simple compact Lie groups, and $C$ is a totally disconnected closed central subgroup, then
$(K_e)_\Cee\cong (T_\Cee\times S_\Cee)/\til{C}$, where $\til{C}$ is a natural image of $C$.  
This is locally compact exactly when both $T_\Cee$ and $S_\Cee$ are locally compact, which
happens exactly what $\rank\what{T}<\infty$ and $S$ is Lie.  Generally, $K_e/Z$ is a quotient
of $S$ by a central subgroup, and $Z_e$ is a quotient of $T$ by a totally disconnected
closed central subgroup so $\rank\what{T}=\rank\what{Z_e}=\rank\what{Z}$, as in Remark
\ref{rem:cotorsion}. 
We observe that the centre of the the space of derivations in the omitted condition (b) is isomorphic to 
$\Ree(\what{T})=\Ree(\what{Z})$.
\end{remark}

\subsection{Some direct sums}
Let for $i$ in an index set $I$, $\Gam_i$ be a non-zero subgroup of $\Que$,
and $\Gam=\bigoplus_{i\in I}\Gam_i$.
Then $\Ree(\Gam)\cong \Ree^{\oplus I}$, so $\rank\Gam=|I|$.  Furthermore,
$\Ree(\Gam)^*\cong\Ree^I$, whose dimension if $2^{|I|}$ is $I$ is infinite; and
$\what{G}^+\cong \Ree^I$, which is locally compact only if $I$ is finite.

 \begin{remark}\label{rem:freab_que}
$\Zee^2$ has a more complicated range of spectra of its Beurling algebras 
 than does $\Que$, i.e.\ $\fL(\Ree^2)$ versus $\fL(\Ree)$.
\end{remark}

Notice that the direct sum structure covers many discretized groups.  For example for 
a non-discrete locally compact field of characteristic zero we have
$\Bbbk_d\cong\Que^{\oplus\mathfrak{c}}$,  whereas if $\Bbbk=\tor(\Bbbk)$
if $\Bbbk$ is non-Archimedean.  Also
$\Tee_d\cong \Que^{\oplus\mathfrak{c}}\oplus \Que/\Zee$,
where $\Que/\Zee\cong\tor(\Tee_d)$.

\subsection{Examples with non-split torsion}\label{ssec:nonsplit}
There are many conditions for a discrete abelian $G$ which
 guarantee that $G\cong G/\tor(G)\oplus\tor(G)$.  For example, this holds if $G$ is either
 finitely generated or divisible.  
 Let us see that this does not always hold.  The following example
 is motivated by one found in some online lecture notes of K. Igusa at Brandeis U.
 
 Let us fix a sequence $(a_k)_{k=1}^\infty$ of natural numbers with
 each at least $2$.  Let $H=\bigcup_{k\in\En}\frac{1}{a_1\dots a_k}\Zee\subseteq\Que$.  We consider
 the quotient of full product discrete groups
 \[
 \Pi=\Zee^\En/P\text{ where }P=\prod_{k\in\En} (a_1\dots a_k)^2\Zee.
\]
Let for $k$ in $\En$
\[
x_k=(0,\dots,0,1,a_{k+1},a_{k+1}a_{k+2},a_{k+1}a_{k+2}a_{k+3},\dots)+P
\]
and $t_k=(0,\dots,0,1,0\dots)+P$, each with $1$ in the $k$th position, and
\[
G=\langle x_k:k\in\En\rangle\subset\Pi
\]
The critical family of relations in this group is
\begin{equation}\label{eq:relation}
x_k-a_{k+1}x_{k+1}=t_k.
\end{equation}
Then $T=\langle t_k:k\in\En\rangle\subseteq\tor(G)$, with each
\[
x_k-(a_{k+1}\dots a_{k+m})x_{k+m}\in T
\]
and hence
\[
\frac{1}{a_1\dots a_k}\mapsto x_k+T\text{ extends to an injection }H\hookrightarrow 
G/T.
\]
This shows that $G/T$ is torsion-free, so $T=\tor(G)$.  Further, $G/T$ 
cannot be isomorphic to a subgroup of $G$, as
no element of $G$ is simultaneously divisible by $a_1\dots a_k$ for all $k$.
The above map is clearly surjective, so
the theory of complexification and weights can then be deduced from Section \ref{ssec:rational}.

It is amusing to compute the dual group.  If $\sig\in\what{G}$, let $z_k=\sig(x_k)$, and 
(\ref{eq:relation}) provides that
\[
z_kz_{k+1}^{-a_{k+1}}=e^{\frac{2\pi i}{(a_1\dots a_k)^2}m_k}\text{ where }
m_k\in\{0,\dots,(a_1\dots a_k)^2-1\}.
\]
Thus
\[
\what{G}=\left\{(z,m)\in\Tee^\En\times \prod_{k\in\En} \Zee/(a_1\dots a_k)^2\Zee:
\text{each }z_k=z_{k+1}^{a_{k+1}}e^{\frac{2\pi i}{(a_1\dots a_k)^2}m_k}\right\}.
\]
Hence we obtain a compact abelian group whose connected component
$G_0=\{z\in\Tee^\En:\text{each }z_k=z_{k+1}^{a_{k+1}}\}\cong\what{H}$
is not a factor group.


\begin{thebibliography}{10}

\bibitem{CartwrightMc} 
D. I. Cartwright, J. R. McMullen, A generalised complexification for compact groups, J. Reine Angew. Math. 331 (1982) 1--15.

\bibitem{dales}
H.G. Dales, Banach Algebras and Automatic Continuity, Clarendon Press, Oxford, 2000.

\bibitem{feichtinger}
H. G. Feichtinger, Gewichtsfunktionen auf lokalkompakten Gruppen,
Sitzber. \"{O}sterr. Akad. Wiss. Abt. II 188 (1979), 451–471 

\bibitem{gharamani_zadeh}
F. Ghahramani and S. Zadeh,
Bipositive isomorphisms between Beurling algebras and between their second dual algebras.
Canad. J. Math. 69 (2017),  3--20.

\bibitem{hewittrI}
E. Hewitt, K. A. Ross, Abstract harmonic analysis. Vol. I.  Springer-Verlag, Berlin-New York, 1979. 

\bibitem{hormander}
L, H\"{o}rmander, Sur la fonction d'appui des ensembles convexes dans un espace localement 
convexe.  Ark. Mat. 3 (1955), 181--186.

\bibitem{hofmannmorris}
K. H Hofmann, S. A. Morris, The  structure of compact groups,  De Gruyter, Berlin, 2006.


\bibitem{kaniuth}
E. Kaniuth, A Course in Commutative Banach Algebras, Springer, New York, 2009.

\bibitem{kuznetsova}
Yu. N. Kuznetsova, Invariant weighted algebras $L_w^p(G)$,
Mat. Zametki 84, (2008) 567--576; translation in
Math. Notes 84 (2008) 529--537.

\bibitem{NariciB}
L. Narici, E. Beckenstein, Topological Vector Spaces, CRC Press., Boca Raton, FL, 2011.

\bibitem{rudin}
W. Rudin, Fourier analysis on groups, John Wiley \& Sons, Inc., New York, 1990.

\bibitem{schaeferwolff}
H. H. Schaefer, M. P. Wolff, Topological vector spaces, Springer-Verlag, New York, 1999.

\end{thebibliography}
\end{document}